\documentclass[12pt,onesided]{article}

\usepackage{graphicx,amssymb,amsmath}
\usepackage{bm}
\usepackage[numbers]{natbib}

\input epsf
\usepackage{theorem}
\usepackage{todonotes}
\usepackage[margin=2cm]{geometry}

\newcommand{\pp}[2]{\frac{\partial #1}{\partial #2}}
\def\eps{\varepsilon}

\def\O{{\mathcal{O}}}

\def\q{{\bf{q}}}
\def\p{{\bf{p}}}
\def\l{{\bf{l}}}
\def\MM#1{{\bf{#1}}}
\def\bphi{{\MM{\phi}}}
\def\Qbar{{\bar{\bf{Q}}}}

\def\u{{\bf{u}}}
\def\v{{\bf{v}}}
\def\w{{\bf{w}}}
\def\Ubar{{\bar{\bf{U}}}}

\def\bzeta{{\bm{\zeta}}}
\DeclareMathOperator{\ad}{ad}
\DeclareMathOperator{\diff}{d}
\newcommand{\difft}{d}
\newcommand{\R}{\mathbb{R}}
\newcommand*{\Id}{\operatorname{Id}}

\title{Stochastic partial differential fluid equations as a diffusive limit of deterministic Lagrangian multi-time dynamics}

\author{C. J. Cotter\thanks{Department of Mathematics, Imperial College, London, UK; {\rm{colin.cotter@imperial.ac.uk}}} \and G.A. Gottwald\thanks{School of Mathematics and Statistics, University of Sydney, Sydney, Australia; {\rm{georg.gottwald@sydney.edu.au}}} \and D. D. Holm\thanks{Department of Mathematics, Imperial College, London, UK; {\rm{d.holm@imperial.ac.uk}}}}

\date{\today}

\begin{document}
\maketitle

%%%%%%%%%%%%%%%%%%%%%%%%%%%%%%%%%%%%%%%%%%%%%%%%%%%%%%%%%%%%

\abstract{In {\em{Holm}, Proc. Roy. Soc. A 471 (2015)} stochastic fluid equations were derived by employing a variational principle with an assumed stochastic Lagrangian particle dynamics. Here we show that the same stochastic Lagrangian dynamics naturally arises in a multi-scale decomposition of the deterministic Lagrangian flow map into a slow large-scale mean and a rapidly fluctuating small scale map. We employ homogenization theory to derive effective slow stochastic particle dynamics for the resolved mean part, thereby justifying stochastic fluid partial equations in the Eulerian formulation. 
%The stochasticity arises through the assumption that the fluctuating part is mildly chaotic and stochastic Lagrangian dynamics can be rigorously derived by means of homogenization. 
To justify the application of rigorous homogenization theory, we assume mildly chaotic fast small-scale dynamics, as well as a centering condition.
The latter requires that the mean of the fluctuating deviations is small, when pulled back to the mean flow.}

\vspace{1cm}
%\noindent
%{\bf{subject}}: mathematical physics, applied mathematics, fluid mechanics\\

\noindent
{\bf{keywords}}: geometric mechanics, stochastic fluid models, stochastic processes, multi-scale fluid dynamics, symmetry reduced variational principles, homogenisation\\

%%%%%%%%%%%%%%%%%%%%%%%%%%%%%%%%%%%%%%%%%%%%%%%%%%%%%%%%%%%%

\section{Introduction}
When studying complex turbulent flows or astrophysical and geophysical fluids, in which physical processes occur over a wide range of spatial and temporal scales, we are faced with the inevitable problem that our limited computational resources will eventually force us to under-represent processes that occur below certain temporal and spatial scales. Processes with spatial scales smaller than the grid-scale and with temporal scales below the applied discrete time step simply  cannot be resolved. However, these unresolved processes may be energetically important and could have a nontrivial impact on the resolved scales. For example, in the study of midlatitude weather, one is generally interested in the dynamics of the synoptic weather systems such as high and low pressure fields with spatial scales of several hundred kilometres and time scales of several days. Numerical simulations must choose a discretisation to resolve those processes. The atmosphere, however, also supports gravity waves, which are the buoyancy oscillations of the stratification surfaces with spatial scales of hundred metres and time scales of several minutes to hours. These unresolved gravity waves play a crucial role in vertical momentum transport; so their omission in a numerical simulation would lead to a misrepresentation of this transport. Similarly, ocean current models do not resolve meso-scale eddies which transport momentum, heat and salinity over large spatial scales. Scientists therefore search for methods to parametrize the influence of these unresolved processes. 

Recently, stochastic partial differential fluid equations have been proposed to model the influence of unresolved scales on the resolved scales of interest \cite{Memin14,Holm15,ResseguierEtAl17a,ResseguierEtAl17b,ResseguierEtAl17c}. These novel approaches introduce stochasticity into the flow map for the Lagrangian particle trajectories, then the noise in the Lagrangian-to-Euler map produces a random Eulerian vector field. This approach results in an interesting form of multiplicative noise in the Eulerian stochastic partial differential fluid equation (SPDE). In the SPDE, the gradient of the noise and its magnitude multiply, respectively, the solution and the gradient of the solution. Thus, the effects of the noise depend on spatial gradients of both the noise and the solution. In \cite{Memin14,ResseguierEtAl17a,ResseguierEtAl17b,ResseguierEtAl17c} the transformation from the Lagrangian to the Eulerian fluid formulation was achieved via a particular version of the Reynold's transport theorem that preserved conservation of energy. 

In \cite{Holm15} a variational principle was employed in the derivation of the SPDE for fluids. The same model was later derived in \cite{CrFlHo2017} from Newton's 2nd Law of Motion and a different version of the Reynold's transport theorem, which included the transformation of the coordinate basis under the Lagrange-to-Euler map. Instead of of preserving the energy, the model derived in \cite{Holm15,CrFlHo2017} has the property of preserving Kelvin's circulation theorem, while introducing a stochastic advecting velocity for the Kelvin loop, which transports the line elements of the Kelvin loop along stochastic Lagrangian trajectories. 

The stochastic Eulerian random field associated with this advecting velocity vector field was prescribed in \cite{Holm15,CrFlHo2017}. Here, we will motivate the stochastic advection of fluid particles by using homogenisation techniques \cite{PavliotisStuart}. In particular, we will show that the same form of stochasticity as in \cite{Holm15,CrFlHo2017} naturally arises, upon splitting the Lagrangian fluid flow map into the composition of mean and fluctuating maps that is standard for generalised Lagrangian mean theories \cite{andrews1978exact,buhler2014waves,holm2002lagrangian,salmon2013alternative}. In particular, when the time scale of the fluctuating map is sufficiently rapid in comparison to the slow mean field dynamics, and when the fast dynamics is sufficiently chaotic, its integrated effect on the slow mean field dynamics is an effective stochastic term, whose variance is determined via the auto-correlation function of the fast fluctuating flow map.\\

The paper is organized as follows. Section~\ref{sec.Holm} briefly summarizes the variational stochastic fluid dynamics introduced in \citet{Holm15}.
%, whose analytical properties were subsequently studied in \cite{CrFlHo2017}. 
Section~\ref{sec.multiscale} then constructs a deterministic multi-scale Lagrangian particle system. In Section~\ref{sec.homo} we employ homogenization to derive the stochastic noise in the Lagrangian dynamics assumed in the variational stochastic fluid approach from the deterministic multi-scale system. We conclude in Section~\ref{sec.summary} with a summary and discussion.

%%%%%%%%%%%%%%%%%%%%%%%%%%%%%%%%%%%%%%%%%%%%%%%%%%%%%%%%%%%%

\section{Variational stochastic fluid equations}
\label{sec.Holm}
\subsection{Variational principles for deterministic fluids}
Governing equations for a variety of ideal fluid dynamics theories can be derived within a variational framework as the Euler-Poincar\'e equations associated with the following variational principle with the action \cite{holm1998euler}
\begin{align}
S(\u,\p,\q) = \int \ell(\u,a) \difft t 
+ \Big\langle \p,\dot{\q}-\u(\q,t) \Big\rangle\,\difft t.
\label{Clebsch-var-princ}
\end{align}
Here, $\q(\l,t)=g(t)\l$ is the Lagrangian trajectory for the time-dependent fluid particle flow map $g(t)$, with $g(0)=\Id$. In this paper, we will denote explicit time dependence with either parentheses, such as $g(t)$, or with subscript $t$, such as as $g_t$; while as usual $\partial_t$ will denote partial time derivative and over-dot the total time derivative. For example, the fluid particle flow map $\q(\l,t)=g(t)\l=g_t\l$ provides the spatial location of the particle with label $\l$ at time $t$.  This is the Lagrange-to-Euler map. Accordingly, the time derivative along the Lagrangian particle trajectory $\q(\l,t)$ satisfies the defining relation for the Eulerian velocity, $\u_t$, given by
\begin{align}
  \label{eq:qdot}
\dot \q(\l,t) = \dot{g}_t\l = (\dot{g}_tg_t^{-1})g_t\l = \u_t\circ g_t\l := \u(\q(\l,t),t)\,.
\end{align}
Thus, in this notation, $\u_t=\dot{g}_tg_t^{-1}$ is an Eulerian time-dependent velocity vector field, $L(g_t,\dot{g}_t,a_0)=L(g_0,\dot{g}_tg_t^{-1},a_0g_t^{-1})=\ell(\u_t,a_t)$ is the reduced Lagrangian, in which $a_t:=a_0\circ g_t^{-1}$ represents a set of advected quantities; such as advected tracers, $\theta = \theta_0\circ g_t^{-1}$, or advected densities, $\rho \diff^3x = (\rho_0\diff^3l)\circ g_t^{-1}$.

The relation \eqref{eq:qdot} is encoded as a constraint in the Clebsch variational formulation in \eqref{Clebsch-var-princ} with the Lagrange
multiplier $\p(\l,t)$, via the inner product
\begin{equation}
  \big\langle \p, \dot{\q}\big\rangle = \int_{\Omega} \p(\l,t)\cdot \dot{\q}(\l,t)
  \diff^n \l,
\end{equation}
where $n$ is the dimension of the domain $\Omega$. The
extremal points of the constrained, reduced action in \eqref{Clebsch-var-princ} then lead to the Euler-Poincar\'e equations
\begin{equation}
\pp{}{t}\frac{\delta \ell}{\delta \u} + \ad^*_\u \frac{\delta
  \ell}{\delta \u} = a\diamond \frac{\delta \ell}{\delta a},
  \label{e.EP}
\end{equation}
where one defines $\ad^*_\u (\v\cdot \diff {\bf x} \otimes \diff^nx) = L_{\u}(\v\cdot \diff {\bf x} \otimes \diff^nx)$, which has components
\begin{equation}
  \ad^*_\u \v = \nabla\cdot(\u\otimes \v) + (\nabla \u)^T\cdot \v
  = L_{\u}\v,
\label{adstar-def-components}  
\end{equation}
and the diamond operator $(\diamond)$ is defined as
\begin{equation}
  \int_\Omega \w \cdot a \diamond \phi \diff^n x = \int_\Omega \phi 
  L_\w \,a\diff^n x,
\end{equation}
for all test functions $\phi$ and $\w$ satisfying the boundary conditions $\w\cdot {\bf n}=0$ on $\partial\Omega$. 

In \eqref{adstar-def-components}, $L_\w a$ denotes the advection
law for the advected quantity in question. Since advected fluid quantities satisfy $a_0=a_t\circ g_t = a(\q(\l,t),t) =:g^*_t a({\bf x}, t)$, where $g^*_t$ denotes the \emph{pullback} of the Lagrange-to-Euler map, we have the advection law
\begin{align}
\begin{split}
\frac{\difft a_0}{\difft t} = 0 & = \frac{\difft }{{\,\difft t}}\big(g^*_t a({\bf x}, t)\big) 
= g^*_t \big((\partial_t   + L_\u ) a_t \big)
= \frac{\difft a_t }{\difft t}
\quad\hbox{along}\quad  \dot \q(\l,t) = \u({\bf x},t)\Big|_{{\bf x}=\q(\l,t)}
\\&= \Big[\partial_t a_t 
+ \dot \q(\l,t)  \cdot  \nabla a_t\Big]_{{\bf x}=\q(\l,t)}  
= \Big[\partial_t a_t + \u({\bf x},t) \cdot  \nabla a_t\Big]_{{\bf x}=\q(\l,t)}
\,.
\end{split}
\label{advect-law}
\end{align}
For example, if $a=\rho \diff^nx$ is a density, then $g^*_t a= g^*_t (\rho \diff^nx) = [\rho({\bf x},t)\diff^nx]_{{\bf x}=\q(\l,t)}$. Thus, the advection law for a density satisfies the continuity equation in coordinates, since in this case the previous formula becomes 
\begin{equation}
0 
= \frac{\difft  }{\difft t} \Big[\rho({\bf x},t)\diff^nx\Big]_{{\bf x}=\q(\l,t)}
= \Big[\big(\partial_t \rho + \u\cdot\nabla\rho 
+ \rho \nabla\cdot\u\big)\diff^nx\Big]_{{\bf x}=\q(\l,t)}.
\end{equation}
Consequently, for a density $\rho \diff^nx$ and a scalar function $\phi$, we find $\rho \diamond \phi = -\rho \nabla \phi $, since
\begin{equation}
  \int_\Omega \w \cdot \rho \diamond \phi \diff^n x 
  = \int_\Omega \phi \,  L_\w \rho \diff^n x
  = \int_\Omega \phi \,  \nabla\cdot(\w \rho) \diff^n x  
  = - \int_\Omega  \w \cdot (\rho \nabla \phi) \, \diff^n x  
  .
\end{equation}
It is slightly more complicated to prove equation \eqref{adstar-def-components}. However, the same rules apply; so
\begin{align}
\begin{split}
g^*_t\big(L_{\u}(\v\cdot \diff {\bf x} \otimes \diff^nx)\big)
&=
\frac{\difft }{{\,\difft t}}\big(g^*_t (\v({\bf x})\cdot \diff {\bf x} \otimes \diff^nx)\big) 
=
\frac{\difft }{{\,\difft t}}\Big[\v({\bf x})\cdot \diff {\bf x} \otimes \diff^nx\Big]_{{\bf x}=\q(\l,t)}
\\& =
\Big[
\big((\u\cdot\nabla) \v + v_j \nabla u^j + \v (\nabla\cdot \u)\big)
\cdot \diff {\bf x} \otimes \diff^nx\Big]_{{\bf x}=\q(\l,t)}
\\&=
\Big[
\Big(\nabla\cdot(\u\otimes \v) + (\nabla \u)^T\cdot \v \Big)
\cdot \diff {\bf x} \otimes \diff^nx\Big]_{{\bf x}=\q(\l,t)}
\\&=
g_t^*\Big[
\Big(\nabla\cdot(\u\otimes \v) + (\nabla \u)^T\cdot \v \Big)
\cdot \diff {\bf x} \otimes \diff^nx\Big]
.
\end{split}
\end{align}
The proof of equation \eqref{adstar-def-components} then follows, because in this case we have the identity,
\begin{align}
L_{\u}(\v\cdot \diff {\bf x} \otimes \diff^nx) 
= {\rm ad}^*_{\u}(\v\cdot \diff {\bf x} \otimes \diff^nx)
\,.
\end{align}
Different choices of $a$ and $\ell$ lead to a huge range of equation
systems in fluid dynamics, including magnetohydrodynamics, geophysical
fluid dynamics and complex fluids, all formulated in terms of their transformation properties under smooth invertible flow maps; see the examples in
\cite{holm1998euler}.  The possibility to eliminate $\p$ and $\q$ from
the equations is a result of the particle relabelling transformation properties of the
action $S$, which leads to Kelvin's circulation theorem, expressed
as \cite{holm1998euler}
\begin{equation}
  \frac{\difft}{\,{{\,\difft t}}}\oint_{C(t)} \frac{1}{\rho}\frac{\delta \ell}{\delta
    \u}\cdot\difft {\bf x} = \oint_{C(t)} \frac{1}{\rho}\left(a\diamond \frac{\delta \ell}{\delta a}\right)\cdot\difft {\bf x},
\label{Kel-Noether-thm}    
\end{equation}
where $C(t)$ is an arbitrary closed loop which is advected by the
velocity $\u$. The proof of \eqref{Kel-Noether-thm} follows immediately by recalling  that the change from Lagrangian to Eulerian variables is the pullback of the Lagrange-to-Euler map, whose time derivative is an advection law. This leads, for example, to Lagrangian conserved potential vorticity, in particular instances of the advected quantity $a$. The equations
also have a conserved energy provided $\ell$ does not depend explicitly on time, which arises from Noether's theorem for  variational principles that are invariant under time translations.

For example, the reduced Lagrangian for the incompressible Euler equations
is 
\begin{equation}
\label{eq:ell euler}
\ell[{\u},\rho] = \int_\Omega \frac12{\rho|{\u}|}^{2} + p(\rho-\rho_0)\diff ^n x,
\end{equation}
where $p$ is a Lagrange multiplier (which will turn out to be the
pressure) enforcing constant density $\rho=\rho_0$. Then, the
Euler-Poincar\'e equations with advected quantities become
\begin{align}
\partial_t(\rho{\u}) + \nabla\cdot({\u}\otimes \rho{\u}) 
+ (\nabla {\u})^T\cdot \rho{\u} 
& = \rho\nabla (-p + |{\u}|/2), \\
\rho_t + \nabla\cdot(\rho{\u}) & = 0, \\
\rho &= \rho_0.
\end{align}
When $\rho_0$ is a constant, these equations simplify to 
\begin{align}
\partial_t {\u} + {\u}\cdot\nabla{\u} = -\nabla p
, \quad\hbox{and}\quad
\nabla\cdot {\u} = 0.
\end{align}

\subsection{Variational principles for stochastic fluids}
In \cite{Holm15,CrFlHo2017}, this framework was extended to include stochastic
parameterisations of unresolved dynamics, by replacing the deterministic
equation \eqref{eq:qdot} with a stochastic equation
\begin{equation}
  \label{eq:dq}
  \difft \MM{q} = \MM{u}(\MM{q},t)\,{{\,\difft t}} + \sum_{i=1}^m \MM{\xi}_i(\MM{q})\circ {\difft W_i},
\end{equation}
%where $\difft$ denotes the stochastic evolution operator, 
where $\MM{\xi}_i(\MM{q})$ are time-independent vector fields, and $W_i$ are
independent stochastic Brownian paths. Inserting the stochastic vector field Ansatz \eqref{eq:dq} into the action $S$
and seeking extrema leads to the stochastic Euler-Poincar\`e equation, and (\ref{e.EP}) is replaced by
\begin{equation}
{\difft}\frac{\delta \ell}{\delta \MM{u}} + \ad^*_{(\MM{u} \,{{\,\difft t}} + \sum_{i=1}^m {\bf
    \MM{\xi}}_i\circ {\difft W_i})}\frac{\delta \ell}{\delta \MM{u}} =
a\diamond \frac{\delta \ell}{\delta a}\,{{\,\difft t}}.
\end{equation}
This equation no longer conserves energy, but it still has a Kelvin circulation
theorem
\begin{equation}
  {\difft} \oint_{C(t)} \frac{1}{\rho}\frac{\delta \ell}{\delta
    \MM{u}}\cdot{\diff} \MM{x} = \oint_{C(t)} a\diamond \frac{\delta \ell}{\delta a}\cdot{\diff} \MM{x} \,{\difft} t,
\end{equation}
where points on the closed curve $C(t)$ now follow the stochastic
vector field in \eqref{eq:dq}. A number of example stochastic geophysical
fluid dynamics models derived in this framework were explored in
\cite{Holm15}.

Returning to our example of the incompressible Euler equation, the 
stochastic Euler-Poincar\'e equation with advected quantities becomes
 \begin{align}
 \begin{split}
 {\difft}(\rho{\MM{u}}) + \nabla\cdot({\difft}\MM{x} \otimes \rho{\MM{u}})% & \\
%\qquad   
+ \left(\nabla {\difft}\MM{x}\right)^T \cdot
 \rho{\MM{u}} & = \rho\nabla (-p + |{\MM{u}}|/2)\,{{\,\difft t}}, \\
 {\difft}\rho + \nabla\cdot(\rho\,{\difft}\MM{x}) & = 0, 
\\
 \rho &= \rho_0,
 \end{split}
 \label{EPcomp-eqn}
 \end{align}
where 
\begin{equation}
{\difft}\MM{x} = {\MM{u}}\,{{\,\difft t}} + \sum_{i=1}^m \MM{\xi}_i\circ {\difft W_i}.
\end{equation}
For constant density, the system \eqref{EPcomp-eqn} simplifies to the stochastic Euler fluid equations analyzed in \cite{CrFlHo2017}; namely,
\begin{align}
{\difft}{\MM{u}} 
+ {\difft}\MM{x}\cdot\nabla{\MM{u}} 
+ \left(\nabla \sum_{i=1}^m \MM{\xi}_i\circ {\difft W_i}\right)^T\!\!\!\! \cdot{\MM{u}}
 = -\nabla p\,{{\,\difft t}}
 , \quad\hbox{with}\quad
\nabla\cdot {\MM{u}} & = 0.
\end{align}

\noindent
{\bf Aim of the paper.}
The aim of this paper is to establish conditions under which the stochastic vector field Ansatz in equation \eqref{eq:dq} may be derived by applying the method of homogenization \cite{Givonetal04,PavliotisStuart}, for the purpose of gaining insight into the situations where such a model can be used in fluid dynamics.

%%%%%%%%%%%%%%%%%%%%%%%%%%%%%%%%%%%%%%%%%%%%%%%%%%%%%%%%%%%%

\section{Stochastic Lagrangian multi-time dynamics}
\label{sec.multiscale}
We postulate a slow-fast evolutionary fluid flow map as the composition of a mean flow $\bar g_t$ depending on slow time $t$ and a rapidly fluctuating flow $g_{t/\eps}^\prime$ associated with the evolution of the fast time scales $t/\eps$, with $\eps\ll1$. We define the flow map associated with the fast scales as the (spatially) smooth invertible map with smooth inverse (i.e., a diffeomorphism, or diffeo for short) given by the sum,
\begin{align}
g_{t/\eps}^\prime = \operatorname{Id} + \,{\zeta}_{t/\eps}
\quad\hbox{where}\quad \eps\ll1
\, .
\label{e.gprime}
\end{align}
The full flow map is taken to be the composition of $\bar g_t$ and $g_{t/\eps}^\prime$, as
\begin{align}
g_t = g_{t/\eps}^\prime\circ \bar g_t = \bar g_t + {\zeta}_{t/\eps}\circ \bar g_t\, .
\label{compmap}
\end{align}
The Lagrangian trajectory of a fluid parcel is then given by $\q(\l,t)=g_t\l$, so that
\begin{align}
\q(\l,t) = \bar{\q}(\l,t) + {\zeta}_{t/\eps}\circ \bar{\q}(\l,t)\,,
\label{e.q}
\end{align}
where the vector $\l$ denotes the fluid label, e.g., the initial condition of the particles. The rapidly fluctuating vector displacement field 
\begin{align}
\bm{\zeta}(\bar{\q}(\l,t),{t/\eps}):={\zeta}_{t/\eps}\circ \bar{\q}(\l,t)
\label{e.vdf}
\end{align}
is defined along the slow large scale resolved trajectory, $\bar\q$. At this point, (\ref{e.q}) may be taken as exact, since it follows directly from the definition of the map ${\zeta}_{t/\eps}$ in (\ref{e.gprime}). 

The tangent to the flow map $g_t$ in \eqref{compmap} at $\q(\l,t)$ along the Lagrangian trajectory \eqref{e.q} defines the Eulerian velocity vector field $\u$, written as
\begin{align}
\dot{g}_t\l = {\q}_t(\l,t) = \u(\q(\l,t),t)\, .
\end{align}
Differentiation of the Lagrangian trajectory (\ref{e.q}) with the assumed fluctuating displacement field (\ref{e.vdf}) yields 
\begin{align}
\u(\bar\q +\zeta_{t/\eps}\circ \bar \q,t) = {\dot{\bar{\q}}}
+  ({\dot{\bar{\q}}}\cdot {\bm{\nabla}_{\bar\q}})\,\bzeta(\bar\q(\l,t),\frac{1}{\eps}t)+ \frac{1}{\eps}\, \partial_t\bzeta \, .
\label{e.qdot}
\end{align}
We shall express the temporal partial derivative $\partial_t\bzeta$ of the fluctuating vector displacement field $\bzeta$ in (\ref{e.vdf}) in terms of its empirical orthogonal eigenfunctions $\bphi_i$, $i=1,2,\dots,M$. We assume that the eigenfunctions $\bphi_i$ are slowly varying in space and are conditioned on the large scale mean flow dynamics, so that we may  write
\begin{align}
\partial_t \bzeta (\bar\q(\l,t),\frac{1}{\eps}t) = 
\sum_{i=1}^{M} \lambda_i(\frac{t}{\eps})\bphi_i(\bar\q(\l,t))
\, .
\label{e.xidot}
\end{align}
Note that the eigenvalues $\lambda_i(t)$ have temporal mean zero. We assume chaotic deterministic  dynamics of the fast fluctuating $\lambda_i(t/\eps)$ with
\begin{align}
\dot \lambda_i = \frac{1}{\eps^2}h_i(\lambda)\, .
\label{e.lambda}
\end{align}
For convenience, let us summarize the deterministic multi-scale system describing the dynamics of the Lagrangian mean position variable $\bar\q$, as follows,
\begin{align}
(\Id + \nabla_{\bar\q}\bzeta)\, {\dot{\bar\q}}&=\u(\bar\q+\zeta_{t/\eps}\circ\bar\q,t)-\frac{1}{\eps}\partial_t\bzeta(\lambda,\bar\q(\l,t))
\label{e.ms1}
\\
\dot\lambda_i&=\frac{1}{\eps^2}h_i(\lambda)
\label{e.ms2}
\\
\partial_t\bzeta&=\sum_{i=1}^{M} \lambda_i(\frac{t}{\eps})\bphi_i(\bar\q(\l,t))
\label{e.ms3}
\end{align}
The next Section employs homogenization theory to show how this set of equations converges for $\eps\to 0$ to solutions of a certain stochastic Lagrangian dynamics on long time scales.

%%%%%%%%%%%%%%%%%%%%%%%%%%%%%%%%%%%%%%%%%%%%%%%%%%%%%%%%%%%%

\section{Diffusive limit of stochastic Lagrangian dynamics}
\label{sec.homo}
Homogenization is a mathematical tool for extracting the statistical behaviour of the slow degrees of freedom in a multi-scale system, by providing an effective stochastic differential equation for the slow degrees of freedom \cite{Givonetal04,PavliotisStuart}. Homogenization represents the integrated effect of the fast (either stochastic, or chaotic) dynamics on the slow variables as noise. Homogenization was initially developed for stochastic multi-scale systems \cite{Khasminsky66,Kurtz73,Papanicolaou76}. However, it has been extended recently to apply for deterministic multi-scale systems with sufficiently chaotic fast dynamics \cite{MelbourneStuart11,GottwaldMelbourne13c,KellyMelbourne17}.
%\footnote{A requirement is that the dynamics can be modelled by a Young tower with summable decay of correlations \cite{MelbourneNicol08,AlvesEtAl11}.} (see \cite{GottwaldMelbourne14} for a discussion).\\ 
In particular, homogenization applies to deterministic multi-scale systems for slow variables $x\in\R^d$ and fast variables $y\in \R^m$ of the form,
\begin{align}
\dot x &=\frac{1}{\eps}f_0(x,y) + f_1(x,y)
\,,\label{e.homo1}
\\
\dot y &= \frac{1}{\eps^2}g_0(y)\, .
\label{e.homo2}
\end{align}
We assume that the vector fields $f_0:\R^d\times\R^m\to\R^d$, $f_1:\R^d\times\R^m\to\R^d$ and $g:\R^m\to\R^m$ satisfy certain regularity conditions. We consider sufficiently chaotic fast $y$-dynamics with a compact chaotic attractor $\Lambda\subset\R^m$ and ergodic invariant probability measure $\mu$. Homogenization requires the so-called {\em{centering condition}} $\int_\Lambda f_0\,d\mu=0$, i.e., that averaging would result in trivial dynamics associated with $f_0$. It is instructive to write the slow dynamics in the integrated form
\begin{align}
x(t) = x(0) + \eps^2 \int_0^{\frac{t}{\eps^2}}f_1(x(\tau),y(\tau)) \difft\tau + \eps\int_0^{\frac{t}{\eps^2}}f_0(x(\tau),y(\tau)) \difft\tau\, ,
\end{align}
where we perform the time integral over the fast time scale $\tau = t/\eps^2$. The middle term involving $f_1$ approaches for $\eps \to 0$ its ergodic mean by means of the law of large numbers. The last term involves summing up sufficiently decorrelated variables $f_0$ with mean zero (the centering condition) and in the limit $\eps \to 0$ leads to Brownian noise by means of the central limit theorem. The randomness comes here from the initial conditions of the fast ergodic deterministic dynamical system. Under the above hypotheses this heuristic argument can be made rigorous \cite{MelbourneStuart11,GottwaldMelbourne13c,KellyMelbourne17}: as $\eps\to 0$ the slow dynamics $x(t)$ converges weakly in the space of continuous functions $C([0,T],\R^d)$ for some $T>0$ to the solution of the homogenized It\^o stochastic differential equation
\begin{align}
\difft X &= F(X) \difft t + \Sigma(X) \difft W_t \;,
\label{e.homog}
\end{align}
where $W_t$ denotes Brownian motion and the drift coefficient is given by
\begin{align}
F(x) = \int f_1(x,y)\,\mu(dy) 
+ \int_0^\infty \difft s \int f_0(x,y)\cdot \nabla f_0(x,y(s)) \,\mu(\diff y) \; ,
\label{e.F}
\end{align}
and the diffusion coefficient is defined by
\begin{align}
\Sigma(X) \Sigma^T(X) &=\int_0^\infty \difft s\int \left( f_0(y) \otimes f_0(y(s)) + f_0(y(s)) \otimes f_0(y) \right)\, \mu(\diff y) \; ,
\label{e.Sigma}
\end{align}
where the outer product between two vectors is defined as $(a\otimes b)_{ij} = a_{i}b_{j}$ (see \cite{PapanicolaouKohler74,IkedaWatanabe,KellyMelbourne17}). 

We remark that, in one spatial dimension $d=1$, the drift term reduces to 
\[
F(X)=\int f_1(x,y)\,\mu(\diff y) + \frac{1}{2} \Sigma(X)\Sigma^\prime(X)\, ,
\]
which implies that the noise is of Stratonovich type. In higher dimensions, the noise is neither of Stratonovich nor of the It\^o type, and the drift term contains nontrivial corrections stemming from $f_0$. The required chaoticity assumption needed to assure the functional central limit theorem, the so called {\em{weak invariance principle}}, 
%\begin{align}
%W^{(\eps)}(t) = \int_0^\infty ds \int f_0(x,y)\cdot \nabla f_0(x,y(s)) \,\mu(dy) 
%\end{align} 
is rather innocuous, and it allows convergence to diffusive laws to be proven for a large class of chaotic dynamical systems. In particular, no assumption about mixing is needed.  For deterministic maps, the convergence to Brownian motion holds when the correlation function of $f_0$ is summable. For flows, it suffices that there is a Poincar\'e map with these properties (irrespective of the mixing properties of the flow). These systems include, but extend far beyond, Axiom A diffeomorphisms and flows, H\'enon-like attractors and Lorenz attractors. Precise statements about their validity can be found in \cite{MelbourneNicol05,MelbourneNicol08,MelbourneNicol09}. We remark that for weakly chaotic dynamics, when the correlations are not summable/integrable, e.g. for systems exhibiting intermittency with long laminar periods, the central limit theorem breaks down and the noise is not Brownian anymore. Instead, the noise is $\alpha$-stable, allowing for unbounded variance and jumps \cite{GottwaldMelbourne13c}.\footnote{We use the terminology strongly and weakly chaotic here in a manner different from the usual distinction between exponential and algebraic decay of correlations; cf. \cite{GottwaldMelbourne13} for further discussion of this difference.}\\

We now show that the deterministic Lagrangian multi-scale system (\ref{e.ms1})-(\ref{e.ms3}) developed in the previous Section is amenable to homogenization and we will derive the associated effective stochastic limit system. For these purposes, we assume sufficiently chaotic dynamics for $\lambda_i$ with ergodic measure $\mu(d\lambda)$. 

First, we show that the centering condition is satisfied. Denoting by angular brackets the average over the fast measure $\mu(d\lambda)$, we express the centering condition as
\begin{align}
\begin{split}
\langle\left(\Id+\nabla\bzeta\right)^{-1}\partial_t\bzeta\rangle 
&= \langle\left(\nabla g_{t/\eps}^\prime\right)^{-1}\partial_t g_{t/\eps}^\prime\rangle \\
&= \left\langle\left(\partial_t \left(\left(g_{t/\eps}^\prime\right)^{-1}\right)\right)\circ g_{t/\eps}^\prime \right\rangle\, = 0.
%= \langle\partial_t\bzeta\rangle +  \O(\eps) = \O(\eps)\, ,
\end{split}
\end{align}
Consequently, the mean displacements in the
fluctuating map vanish, when pulled back to (\emph{i.e.}, transformed back to vectors relative to) the mean
Eulerian coordinates. In fact, this condition simply defines how we
take the mean; namely, it is defined so that  this condition is satisfied.

Homogenization theory for deterministic multi-scale systems as developed in \cite{MelbourneStuart11,GottwaldMelbourne13c,KellyMelbourne17} assures that the slow dynamics of the deterministic multi-scale Lagrangian particle system (\ref{e.ms1})-(\ref{e.ms3}) is on long time-scales of $\O(1/\eps^2)$ described by the stochastic differential equation
\begin{align}
\difft\Qbar = \Ubar(\Qbar) {\,\difft t} + \MM{\sigma}(\Qbar) \difft W_t\, ,
\label{e.Qbar}
\end{align}
where the drift term is given by
\begin{align}
\label{e.drift}
\Ubar(\Qbar)&= \langle \left(\Id+\nabla \bzeta\right)^{-1}u(\Qbar(\l,t)+\bzeta,t)\rangle
+\langle \partial_t\bzeta \nabla_{\Qbar}\partial_t \bzeta\rangle 
%\nonumber
%\\
%&= \langle u(\bar\q(\l,t)+\bzeta,t)\rangle + \O(\eps)\, ,
\end{align}
%\todo[inline]{$\langle \partial{\bf{\zeta}} \nabla_{\Qbar} \partial{\bf{\zeta}}\rangle$ term is ${\mathcal{O}(\eps)$?}
and the diffusion tensor is given by
\begin{align}
\frac{1}{2}\sigma\sigma^T &= \int_0^\infty \difft s\,
\sum_{i,j=1}^M
\langle
\lambda_i(0)\lambda_j(s)
\left(\Id+\nabla \bzeta\right)^{-1}P_{ij}
\left(\Id+\nabla \bzeta\right)^{-T}
\rangle \,,
%\nonumber
%\\
%&= \int_0^\infty ds\,
%\Phi(\bar\q)
%\langle
%\lambda(0)\lambda^T(s)
%\rangle
%\Phi(\bar\q)^T + \O(\eps) \, ,
\label{e.diffusion}
\end{align}
where $P_{ij}=\bphi_i(\bar\q(\l,t))\bphi_j(\bar\q(\l,t))^T$, and $\difft\MM{W}_t=(\difft W_1,\ldots,\difft W_M)$ is a vector of independent Brownian motions (cf. (\ref{e.F}) and (\ref{e.Sigma})). Upon defining $\MM{\xi}_i$ as the $i$th row of $\sigma$, this result matches the Ansatz \eqref{eq:dq} introduced in the stochastic variational framework \cite{Holm15,CrFlHo2017}. 

Note that because of the multiplicative nature of the noise, the drift (\ref{e.drift}) is not just given by the average mean flow $ \langle \left(\Id+\nabla \bzeta\right)^{-1}u(\Qbar(\l,t)+\bzeta,t)\rangle$; there is an additional drift component provided by the fluctuating vector displacement field.

%%%%%%%%%%%%%%%%%%%%%%%%%%%%%%%%%%%%%%%%%%%%%%%%%%%%%%%%%%%%

\section{Summary and discussion}
\label{sec.summary}
The diffusive limit we have described relies on homogenization theory. In order to apply it we had to make several assumptions. First, we made an \emph{envelope assumption}, that the fluctuating vector displacement field $\bzeta$  varies rapidly in time, but only slowly in space. Second, we made a \emph{chaoticity assumption}, whereby the temporal derivative of $\bzeta$ evolves according to some unspecified chaotic dynamics. The microscopic details of the fast chaotic dynamics are not relevant for the emergence of stochastic Lagrangian particle dynamics, but they do feature in the diffusion tensor in equation (\ref{e.diffusion}). In order to build reliable stochastic coarse-grained fluid models, one needs to determine the drift and diffusion terms (\ref{e.drift}) and (\ref{e.diffusion}), respectively. The advantage of the approach taken in \cite{Holm15,CrFlHo2017} and in this paper is that the fluctuating small-scale field $\bzeta$ is Eulerian. One can therefore obtain the statistics by investigating time-series of Eulerian flow fields in terms of their empirical orthogonal eigenfunctions $\bphi$ and their temporal evolution. The data may be generated either by analysis of  numerical high-resolution simulations, or from observational data.\\

The diffusive limit in the Lagrangian dynamics for the slowly varying
mean flow then forms the basis for the variational approach proposed
in \cite{Holm15,CrFlHo2017} to derive stochastic partial differential equations
for fluid systems. The same procedure can, however, also be applied to
give a deterministic multi-scale backbone for the stochastic fluid
equations derived in
\cite{Memin14,ResseguierEtAl17a,ResseguierEtAl17b,ResseguierEtAl17c}.

Equation \eqref{e.drift} shows that there are two contributions to
the mean velocity $\bar{{\u}}$. The first term arises from averaging the full
velocity over the fluctuations, and the second term arises purely from
the structure of the statistics of the time-derivative of the
fluctuation vector $\MM{\zeta}$. The second term shows that the
fluctuating dynamics affects the mean flow. 

The present 
derivation shows that the stochastic fluid equations can be regarded 
as modelling the evolution mean velocity $\bar{{\u}}$ for a fluid that
has a fluctuating velocity that is $\mathcal{O}(1/\epsilon)$. These equations
also provide a possible way to calculate the EOF eigenvectors $\xi_i$ describing spatial correlations  for the stochastic model proposed in \cite{Holm15,CrFlHo2017}. For this purpose, one should compute a large number of Lagrangian trajectories, 
then apply a low-pass spatial filter to them, and compute time series of displacement vectors ${\bf \zeta_t}$ for each trajectory. The values of ${\bf \zeta_t }$ may be binned according coarse-grained grid boxes, and then combined in the formula for $\sigma$.
One can then compute EOFs of $\sigma$, which, in turn, give the EOF eigenvectors $\xi_i$ for the stochastic model proposed in \cite{Holm15,CrFlHo2017}.
The performance of the stochastic model constructed this way is currently
being evaluated, using data assimilation algorithms.

%%%%%%%%%%%%%%%%%%%%%%%%%%%%%%%%%%%%%%%%%%%%%%%%%%%%%%%%%%%%

%\funding{Insert funding text here.}

\section*{Acknowledgments}
This work was partially supported by the EPSRC Standard Grant EP/N023781/1. We would like to thank Greg Pavliotis stimulating discussions. GAG acknowledges the hospitality of Imperial College where this work started.

%%%%%%%%%%%%%%%%%%%%%%%%%%%%%%%%%%%%%%%%%%%%%%%%%%%%%%%%%%%%

%\section*{References}
\bibliographystyle{natbib}
\bibliography{bibliography.bib}

\begin{thebibliography}{}

\bibitem[Andrews and McIntyre(1978)Andrews and McIntyre]{andrews1978exact}
Andrews, D.~G. and McIntyre, M. (1978).
\newblock An exact theory of nonlinear waves on a {L}agrangian-mean flow.
\newblock {\em Journal of Fluid Mechanics\/}, {\bf 89}(04), 609--646.

\bibitem[B{\"u}hler(2014)B{\"u}hler]{buhler2014waves}
B{\"u}hler, O. (2014).
\newblock {\em Waves and mean flows\/}.
\newblock Cambridge University Press.

\bibitem[Crisan {\em et~al.}(2017)Crisan, Flandoli, and Holm]{CrFlHo2017}
Crisan, D.~O., Flandoli, F., and Holm, D.~D. (2017).
\newblock Solution properties of a 3{D} stochastic {E}uler fluid equation.
\newblock {\em arXiv preprint arXiv:1704.06989\/}.

\bibitem[Givon {\em et~al.}(2004)Givon, Kupferman, and Stuart]{Givonetal04}
Givon, D., Kupferman, R., and Stuart, A. (2004).
\newblock Extracting macroscopic dynamics: Model problems and algorithms.
\newblock {\em Nonlinearity\/}, {\bf 17}(6), R55--127.

\bibitem[Gottwald and Melbourne(2013a)Gottwald and
  Melbourne]{GottwaldMelbourne13}
Gottwald, G.~A. and Melbourne, I. (2013a).
\newblock {A Huygens principle for diffusion and anomalous diffusion in
  spatially extended systems}.
\newblock {\em Proc. Natl. Acad. Sci. USA\/}, {\bf 110}, 8411--8416.

\bibitem[Gottwald and Melbourne(2013b)Gottwald and
  Melbourne]{GottwaldMelbourne13c}
Gottwald, G.~A. and Melbourne, I. (2013b).
\newblock Homogenization for deterministic maps and multiplicative noise.
\newblock {\em Proceedings of the Royal Society A: Mathematical, Physical and
  Engineering Science\/}, {\bf 469}(2156).

\bibitem[Holm(2002)Holm]{holm2002lagrangian}
Holm, D.~D. (2002).
\newblock Lagrangian averages, averaged {L}agrangians, and the mean effects of
  fluctuations in fluid dynamics.
\newblock {\em Chaos: An Interdisciplinary Journal of Nonlinear Science\/},
  {\bf 12}(2), 518--530.

\bibitem[Holm(2015)Holm]{Holm15}
Holm, D.~D. (2015).
\newblock Variational principles for stochastic fluid dynamics.
\newblock {\em Proc. Roy. Soc. A.}, {\bf 471}(2176), 20140963, 19.

\bibitem[Holm {\em et~al.}(1998)Holm, Marsden, and Ratiu]{holm1998euler}
Holm, D.~D., Marsden, J.~E., and Ratiu, T.~S. (1998).
\newblock The {E}uler--{P}oincar{\'e} equations and semidirect products with
  applications to continuum theories.
\newblock {\em Advances in Mathematics\/}, {\bf 137}(1), 1--81.

\bibitem[Ikeda and Watanabe(1981)Ikeda and Watanabe]{IkedaWatanabe}
Ikeda, N. and Watanabe, S. (1981).
\newblock {\em {Stochastic Differential Equations and Diffusion Processes}\/},
  volume~24 of {\em North-Holland Mathematical Library\/}.
\newblock North-Holland, New York.

\bibitem[Kelly and Melbourne(2017)Kelly and Melbourne]{KellyMelbourne17}
Kelly, D. and Melbourne, I. (2017).
\newblock Deterministic homogenization for fast--slow systems with chaotic
  noise.
\newblock {\em Journal of Functional Analysis\/}, {\bf 272}(10), 4063 -- 4102.

\bibitem[Khasminsky(1966)Khasminsky]{Khasminsky66}
Khasminsky, R.~Z. (1966).
\newblock {On stochastic processes defined by differential equations with a
  small parameter}.
\newblock {\em Theory of Probability and its Applications\/}, {\bf 11},
  211--228.

\bibitem[Kurtz(1973)Kurtz]{Kurtz73}
Kurtz, T.~G. (1973).
\newblock A limit theorem for perturbed operator semigroups with applications
  to random evolutions.
\newblock {\em Journal of Functional Analysis\/}, {\bf 12}(1), 55--67.

\bibitem[Melbourne and Nicol(2005)Melbourne and Nicol]{MelbourneNicol05}
Melbourne, I. and Nicol, M. (2005).
\newblock Almost sure invariance principle for nonuniformly hyperbolic systems.
\newblock {\em Commun. Math. Phys.}, {\bf 260}, 131--146.

\bibitem[Melbourne and Nicol(2008)Melbourne and Nicol]{MelbourneNicol08}
Melbourne, I. and Nicol, M. (2008).
\newblock Large deviations for nonuniformly hyperbolic systems.
\newblock {\em Trans. Amer. Math. Soc.}, {\bf 360}(12), 6661--6676.

\bibitem[Melbourne and Nicol(2009)Melbourne and Nicol]{MelbourneNicol09}
Melbourne, I. and Nicol, M. (2009).
\newblock A vector-valued almost sure invariance principle for hyperbolic
  dynamical systems.
\newblock {\em Annals of Probability\/}, {\bf 37}, 478--505.

\bibitem[Melbourne and Stuart(2011)Melbourne and Stuart]{MelbourneStuart11}
Melbourne, I. and Stuart, A. (2011).
\newblock A note on diffusion limits of chaotic skew-product flows.
\newblock {\em Nonlinearity\/}, {\bf 24}, 1361--1367.

\bibitem[M\'emin(2014)M\'emin]{Memin14}
M\'emin, E. (2014).
\newblock Fluid flow dynamics under location uncertainty.
\newblock {\em Geophys. Astrophys. Fluid Dyn.}, {\bf 108}(2), 119--146.

\bibitem[Papanicolaou(1976)Papanicolaou]{Papanicolaou76}
Papanicolaou, G.~C. (1976).
\newblock Some probabilistic problems and methods in singular perturbations.
\newblock {\em Rocky Mountain Journal of Mathematics\/}, {\bf 6}(4), 653--674.

\bibitem[Papanicolaou and Kohler(1974)Papanicolaou and
  Kohler]{PapanicolaouKohler74}
Papanicolaou, G.~C. and Kohler, W. (1974).
\newblock Asymptotic theory of mixing stochastic ordinary differential
  equations.
\newblock {\em Comm. Pure Appl. Math.}, {\bf 27}, 641--668.

\bibitem[Pavliotis and Stuart(2008)Pavliotis and Stuart]{PavliotisStuart}
Pavliotis, G.~A. and Stuart, A.~M. (2008).
\newblock {\em {Multiscale Methods: Averaging and Homogenization}\/}.
\newblock Springer, New York.

\bibitem[Resseguier {\em et~al.}(2017a)Resseguier, M\'emin, and
  Chapron]{ResseguierEtAl17a}
Resseguier, E., M\'emin, E., and Chapron, B. (2017a).
\newblock Geophysical flows under location uncertainty, {P}art {I}: {R}andom
  transport and general models.
\newblock {\em Geophys. Astrophys. Fluid Dyn.}, {\bf 111}(3), 149--176.

\bibitem[Resseguier {\em et~al.}(2017b)Resseguier, M\'emin, and
  Chapron]{ResseguierEtAl17b}
Resseguier, E., M\'emin, E., and Chapron, B. (2017b).
\newblock Geophysical flows under location uncertainty, {P}art {II}:
  {Q}uasigeostrophic models and efficient ensemble spreading.
\newblock {\em Geophys. Astrophys. Fluid Dyn.}, {\bf 111}(3), 177--208.

\bibitem[Resseguier {\em et~al.}(2017c)Resseguier, M\'emin, and
  Chapron]{ResseguierEtAl17c}
Resseguier, E., M\'emin, E., and Chapron, B. (2017c).
\newblock Geophysical flows under location uncertainty, {P}art {III}: {SQG} and
  frontal dynamics under strong turbulence.
\newblock {\em Geophys. Astrophys. Fluid Dyn.}, {\bf 111}(3), 209--227.

\bibitem[Salmon(2013)Salmon]{salmon2013alternative}
Salmon, R. (2013).
\newblock An alternative view of generalized {L}agrangian mean theory.
\newblock {\em Journal of Fluid Mechanics\/}, {\bf 719}, 165--182.

\end{thebibliography}

%\bibliographystyle{agsm}
%\bibliography{bibliography}

%%%%%%%%%%%%%%%%%%%%%%%%%%%%%%%%%%%%%%%%%%%%%%%%%%%%%%%%%%%%

\end{document}